\documentclass[11pt]{article}
\usepackage{amsmath,amssymb,latexsym}
\newtheorem{theo}{Theorem}
\newtheorem{coro}{Corollary}

\newtheorem{lemm}{Lemma}
\newtheorem{prop}{Proposition}

\def\iI{\stackrel{\circ}{I}}
\newcommand{\CQFD}
{%
\mbox{}%
\nolinebreak%
\hfill%
\rule{2mm}{2mm}%
\medbreak%
\par%
}

\makeatletter\@addtoreset{equation}{section}\makeatother
\begin{document}
\title{A new extension of bivariate FGM copulas}
\author{C\'ecile Amblard$^{1}$ \& St\'ephane Girard$^{2}$
}
\date{
\small
$^{1}$ Institut de l'Ing\'enierie de l'Information de Sant\'e, \\
TIMC - TIMB, Facult\'e de M\'edecine, Pav. D, \\
38706 La Tronche Cedex, France. \\
E-mail : {\tt Cecile.Amblard@imag.fr}\\
$^{2}$ INRIA Rh\^one-Alpes, team Mistis,\\
655, avenue de l'Europe, Montbonnot,\\
38334 Saint-Ismier Cedex, France. \\ 
 E-mail: {\tt Stephane.Girard@inrialpes.fr}
}
\maketitle

\begin{abstract}
\noindent We propose a new family of copulas generalizing the
 Farlie-Gumbel-Morgenstern family and generated by two univariate functions.
The main feature of this family is to permit the modeling of high positive
dependence. In particular, it is established that the range of the
Spearman's Rho is $[-3/4,1]$ and that the upper tail dependence coefficient can reach any value in $[0,1]$. Necessary and sufficient conditions
are given on the generating functions in order to obtain
various dependence properties.
Some examples of parametric subfamilies are provided.
\\



\end{abstract}

                   \section{Introduction}

\noindent A bivariate copula defined on the unit square $I^2:=[0,1]^2$ is a bivariate cumulative distribution function 
(cdf) with univariate uniform $I$ margins {\it i.e.} verifying the 
following three properties:
\begin{itemize}
\item [(P1)] $C(u,0)=C(0,v)=0$, $\forall (u,v)\in I^2$,
\item [(P2)] $C(u,1)=u$ and $C(1,v)=v$, $\forall (u,v)\in I^2$,
\item [(P3)] $\Delta(u_1, u_2, v_1, v_2):=C(u_2,v_2)-C(u_2,v_1)-C(u_1,v_2)+C(u_1,v_1)\geq 0$, $\forall (u_1, u_2, v_1, v_2)\in I^4$, such that $u_1\leq u_2$ and $v_1\leq v_2$.
\end{itemize}
Sklar's Theorem~\cite{SKLAR} states that any
bivariate distribution with cdf $H$
and marginal cdf $F$ and $G$ can be written
$H(x,y)=C(F(x),G(y))$, where $C$ is a copula. This result justifies the use of copulas for building bivariate distributions.

One of the most popular parametric family of copulas is the Farlie-Gumbel-Morgenstern
(FGM) family defined when $\theta\in[-1,1]$ by
\begin{equation}
\label{eqFGM}
C_\theta^{\mbox{\tiny FGM}}(u,v)=uv+\theta u(1-u)v(1-v),
\end{equation}
and studied in~\cite{FAR,GUMBEL,MORG}. A well-known limitation to this
family is that it does not allow the modeling of large dependences
since Spearman's Rho is
limited to $\rho\in[-1/3,1/3]$. Basing on this remark, more general
copulas have been introduced in 1960 by Sarmanov~\cite{SAR},
\begin{equation}
\label{eqSarmanov}
C_{\theta,\phi,\psi}^{\mbox{\tiny Sarmanov}}(u,v)=uv+\theta \phi(u)\psi(v),
\end{equation}
an re-discovered in 2004 by Rodr\'{\i}guez-Lallena and \'Ubeda-Flores~\cite{LLALENA2}.
See~\cite{KIM} for an extension of this model.
Properties of Sarmanov copulas are studied in~\cite{LEE,SHU}.
Unfortunately, characterization of admissible parameters $\theta$
and functions $\phi$ and $\psi$ is not tractable to obtain
closed-form bounds on the corresponding Spearman's Rho.
Thus, several parametric sub-families of~(\ref{eqSarmanov}) were
introduced. 
In~\cite{NELSEN97}, it is remarked that copulas with quadratic sections~\cite{QUESADA}
are not able to model large dependences. 
Copulas with cubic sections are thus introduced, with the conclusion that
copulas with higher order polynomial sections would increase the dependence degrees but simultaneously the complexity of the model. 
In~\cite{HUE}, two kernel extensions of FGM copulas are studied
\begin{equation}
\label{eqHK1}
C_{\theta,\gamma}^{\mbox{\tiny HK1}}(u,v)=uv+\theta u(1-u)^\gamma v(1-v)^\gamma,
\end{equation}
for $\gamma\geq 1$ and
\begin{equation}
\label{eqHK2}
C_{\theta,\gamma}^{\mbox{\tiny HK2}}(u,v)=uv+\theta u(1-u^\gamma) v(1-v^\gamma),
\end{equation}
for $\gamma\geq 1/2$.
It is shown that Spearman's Rho can be
increased up to approximatively $0.39$ while the lower bound remains
$-{1/3}$. Another similar extension is
\begin{equation}
\label{eqLX}
C_{\theta,p,q}^{\mbox{\tiny LX}}(u,v)=uv+\theta u^p(1-u)^q v^p(1-v)^q,
\end{equation}
see~\cite{LAI}.
Copulas~(\ref{eqHK1}) and~(\ref{eqHK2})
are particular cases of Bairamov-Kotz family~\cite{BK} defined by
\begin{equation}
\label{eqBK}
C_{\theta,p,q,n}^{\mbox{\tiny BK}}(u,v)=u^pv^p\left[1+\theta (1-u^q)^n (1-v^q)^n\right],
\end{equation}
and with associated Spearman's Rho $\rho\in[-0.48, 0.501594]$.
Moreover, it has been remarked in~\cite{HUE2} that dependence degrees arbitrarily close to
$\pm 1$ cannot be obtained with polynomial functions of fixed degree.
An alternative approach to generalize the FGM family of copulas is to consider
the semi-parametric family of symmetric copulas defined by
\begin{equation}
\label{eqcopunous}
C_{\theta,\phi}^{\mbox{\tiny SP}}(u,v)=uv+\theta \phi(u)\phi(v),
\end{equation}
with $\theta\in[-1,1]$. It was first introduced in~\cite{LLALENA},
and extensively studied in~\cite{Nous,Nous2}. 
Let us precise that, in this paper and in accordance with~\cite{NELSEN06},
page 38, a copula $C$ is said to be symmetric if $C(u,v)=C(v,u)$
for all $(u,v)\in I^2$.
Clearly, this family also includes the FGM copulas~(\ref{eqFGM})
(which contains all copulas with both horizontal and vertical quadratic 
sections~\cite{QUESADA}), the parametric family 
of symmetric copulas with cubic sections proposed in~\cite{NELSEN97},
equation~(4.4), and kernel families~(\ref{eqHK1}),~(\ref{eqHK2}) and~(\ref{eqLX}).
It can be shown that, for a properly chosen function $\phi$,
the range of Spearman's Rho is extended to $[-3/4,3/4]$,
whereas the upper tail dependence coefficient is always null.
We refer to~\cite{FISC} for a very interesting method for
constructing admissible functions $\phi$.

In this paper, we propose an extension of the $C_{\theta,\phi}^{\mbox{\tiny SP}}$
family where $\theta$ is a univariate function. This modification
allows the introduction of a singular component concentrated on
the diagonal $v=u$. Consequently, the
modeling of strong positive dependences is possible since 
this new family can take into account the extremal case of positive functional
dependence between two random variables.
Moreover, arbitrary upper tail dependence
coefficients can be reached in $[0,1]$.
The new family is described in Section~\ref{secdef} and the associated
Spearman's Rho and tail dependence coefficients are studied in
Section~\ref{measure}. Section~\ref{dependence} is dedicated to the
dependence properties of this new family of copulas. Finally, some
examples of copulas taken in this family are provided in Section~\ref{secsub}. 
Lemmas are postponed to the appendix.

            \section{Definition and basic properties}
\label{secdef}

We consider the family of functions defined on $I^2$ by:
\begin{equation}
\label{defcop}
C_{\theta,\phi}(u,v)=uv+\theta(\max(u,v)) \phi(u)\phi(v),
\end{equation}
where $\phi$ and $\theta$ are two $I\to {\mathbb R}$ 
continuously differentiable functions. 
Remark that, if $\theta$ or $\phi$ is the null function on $I$
 then $C_{\theta,\phi}=\Pi$, the independent copula.  
In the sequel, we thus assume that
$\phi$ vanishes at most on isolated points of $I$,
and that $\theta$ is not the zero function on $I$.
The next theorem gives sufficient and necessary conditions on $\phi$ and $\theta$ to ensure that $C_{\theta,\phi}$ is a copula. 
\begin{theo}
\label{theocop}
$C_{\theta,\phi}$ is a copula if and only if $\phi$ and $\theta$ satisfy the following conditions~:
\begin{itemize}
\item [(a)] $\phi(0)=0$,
\item [(b)] $\phi(1)\theta(1)=0$,
\item [(c)] $\phi'(u)(\theta\phi)'(v)\geq -1$ for all $0<u\leq v<1$,
\item [(d)] $\theta'(u)\leq 0$ for all $u \in \iI=(0,1)$.
\end{itemize}
\end{theo}
\noindent {\bf Proof: } 
The proof involves four steps.\\
1. First, it is clear that (P1) $\Leftrightarrow$
(a)  and (P2) $\Leftrightarrow$ (b).\\
2. Second, we show that (P3) $\Rightarrow$ (c). 
To this end, consider $0<u_1<u_2\leq v_1<v_2<1$. In this case,
$\Delta(u_1, u_2, v_1, v_2)$ can be rewritten as
\begin{equation}
\label{eqdelta}
\Delta(u_1, u_2, v_1, v_2)=(u_2-u_1)(v_2-v_1) +(\phi(u_2)-\phi(u_1))(\theta(v_2)\phi(v_2)- \theta(v_1)\phi(v_1)),
\end{equation}
and thus  $\Delta(u_1, u_2, v_1, v_2) \geq 0$ implies
$$
\frac{\phi(u_2)-\phi(u_1)}{u_2-u_1}\frac{\theta(v_2)\phi(v_2)
-\theta(v_1)\phi(v_1)}{v_2-v_1}\geq -1.
$$
Letting $u_1\to u_2^-$ 
and $v_1\to v_2^-$ 
in the previous inequality yields (c).\\
3. Similarly, we now show that (P3) $\Rightarrow$ (d). 
Taking $0<u<v<1$, we have
\begin{eqnarray*}
\Delta(u,v,u,v)&=&(v-u)^2+\theta(v)\phi^2(v)+\theta(u)\phi^2(u)-2\phi(u)\phi(v)\theta(v),\\
                &=&(v-u)\left[(v-u)+\theta(v)\phi(v)\frac{\phi(v)-\phi(u)}{v-u}-\phi(u)
\frac{(\theta\phi)(v)-(\theta\phi)(u)}{v-u}\right].
\end{eqnarray*}
Letting $u\to v^-$ in the inequality
$\Delta(u,v,u,v)\geq 0$ yields
$$
\theta(v)\phi(v)\phi'(v)-\phi(v)(\theta\phi)'(v)\geq 0,
$$
which is equivalent to $\phi^2(v)\theta'(v)\leq 0$.
This implies that $\theta'(v)\leq 0$ for all $v\in I$ such that $\phi(v)\neq 0$.
Since $\phi$ vanishes only on isolated points and $\theta'$ is continuous,
(d) is proved.\\
4. Finally, it remains to prove that (c, d) $\Rightarrow$ (P3).
Let $(u_1,u_2,v_1,v_2)\in I^4$ such that $u_1\leq u_2$ and $v_1\leq v_2$.
Let us denote by $R$ the rectangle $[u_1,u_2]\times[v_1,v_2]$,
by $T$ the triangle with vertices $(0,0)$, $(0,1)$ and $(1,1)$,
and by $\bar{T}$ the triangle $\bar{T}=I^2\setminus T$.
Suppose (c, d) hold and let us prove that $\Delta(R):=\Delta(u_1,u_2,v_1,v_2)\geq 0$.
Three cases are considered. 
\begin{enumerate}
\item[(i)] If $R\subset T$, {\it i.e.}  $u_2\leq v_1$ then $\Delta(u_1, u_2, v_1, v_2)$
can be written as in~(\ref{eqdelta}) and
the mean value theorem entails that there exist $u\in(u_1,u_2)$ and $v\in(v_1,v_2)$ such that
\begin{eqnarray*}
\Delta(u_1, u_2, v_1, v_2)&=&(u_2-u_1)(v_2-v_1)\left[ 1  +
\frac{\phi(u_2)-\phi(u_1)}{u_2-u_1} \frac{\theta(v_2)\phi(v_2)-\theta(v_1)\phi(v_1)}{v_2-v_1}
\right]\\
&=& (u_2-u_1)(v_2-v_1)\left[ 1  +
\phi'(u) (\theta\phi)'(v) \right]\geq 0,
\end{eqnarray*}
as a consequence of (c).
\item[(ii)] If $R\subset \bar{T}$, then symmetry considerations
show that $\Delta(u_1,u_2,v_1,v_2)=\Delta(v_1,v_2,u_1,u_2)\geq 0$ from the case (i).
\item[(iii)] If $R\cap T\neq \emptyset$ and $R\cap \bar{T}\neq \emptyset$, then $R$ can be decomposed as non-overlapping 
rectangles $R=R_1\cup R_2\cup R_3$
such that $R_1\subset T$ or $R_1\subset \bar{T}$, $R_2\subset T$ or $R_2\subset \bar{T}$ and $R_3$ is
a square of the form $[u,v]\times[u,v]$. Thus,
$\Delta(R)=\Delta(R_1)+\Delta(R_2)+\Delta(R_3)$ and 
(i), (ii) entail that $\Delta(R_1)\geq 0$ and $\Delta(R_2)\geq 0$.
Let us focus on $\Delta(R_3)$:
\begin{eqnarray*}
\Delta(u,v,u,v)&=&(v-u)^2+\theta(v)\phi(v)[\phi(v)-\phi(u)]-\phi(u)[\theta(v)\phi(v)-\theta(u)\phi(u)]\\
		&=&(v-u)^2+\int_u^v[\theta(v)\phi(v)\phi'(t)-\phi(u)(\theta\phi)'(t)]dt.
\end{eqnarray*}
Note that (c) implies that for all $0<t<v<1$, 
$$
\int_t^v \phi'(t) (\theta\phi)'(y)dy\geq \int_t^v -1 dy = t-v,
$$
and thus $\phi'(t)(\theta\phi)(v)\geq t-v+\theta(t)\phi(t)\phi'(t).$
Similarly, (c) shows that for all $0<u<t<1$,
$$
(\theta\phi)'(t)\displaystyle\int_u^t\phi'(x)dx\geq u-t,
$$ 
and thus $ -(\theta\phi)'(t)\phi(u)\geq u-t-(\theta\phi)'(t)\phi(t)$.
It follows that
$$
\Delta(u,v,u,v)\geq\int_u^v[(\theta\phi)(t)\phi'(t)-\phi(t)(\theta\phi)'(t)]dt
=-\int_u^v\theta'(t)\phi^2(t)dt\geq 0
$$
under condition (d). As a conclusion, $\Delta(R)\geq 0$ and (P3) is proved.
\CQFD
\end{enumerate}

\noindent Note that (b) is true if $\phi(1)=0$ 
or $\theta(1)=0$. We refer to Section~\ref{secsub} for a detailed
study of the corresponding sub-families.
Although the copula $C_{\theta,\phi}$ has full support $I^2$,
the following proposition shows that, in general, it is neither
absolutely continuous, nor singular.  
\begin{prop}
\label{propsing}
The copula $C_{\theta,\phi}$ has both absolutely continuous 
and singular components $A_{\theta,\phi}$ and $S_{\theta,\phi}$,
respectively, given by
$$
A_{\theta,\phi}(u,v)=uv+\theta(\max(u,v))\phi(u)\phi(v)+\int_{0}^{\min(u,v)} (\theta'\phi^2)(t)dt,
$$
and
$$
S_{\theta,\phi}(u,v)=-\int_{0}^{\min(u,v)} (\theta'\phi^2)(t)dt.
$$
\end{prop}
\noindent {\bf Proof: } The absolutely continuous component of
 $C_{\theta,\phi}$ is given by
 $$
 A_{\theta,\phi}(u,v)=\int_0^u\int_0^v 
 \frac{\partial^2} {\partial s\partial t} C_{\theta,\phi}(s,t)dt ds,
$$
with 
$$
 \frac{\partial^2} {\partial s\partial t} C_{\theta,\phi}(s,t)
 = 1 + (\theta\phi)'(\max(s,t)) \phi'(\min(s,t)).
$$
Assume for instance $v\geq u$. Then, $A_{\theta,\phi}$ can be written as
\begin{eqnarray*}
 A_{\theta,\phi}(u,v)&=& uv + \int_0^u\int_0^s (\theta\phi)'(s)\phi'(t)dt ds
 + \int_0^u\int_s^v (\theta\phi)'(t)\phi'(s)dt ds\\
 &=& uv + (\theta\phi)(v)\phi(u) + \int_0^u (\theta'\phi^2)(s)ds\\
 &=& C_{\theta,\phi}(u,v) + \int_0^u (\theta'\phi^2)(s)ds,
\end{eqnarray*}
and the conclusion follows.  The case $v<u$ is similar.  
\CQFD
\noindent Thus, the mass of the singular component is concentrated on the 
first diagonal of the square $I^2$.  Denoting by $(U,V)$ a uniform
random pair on $I^2$ with copula $C_{\theta,\phi}$, we have
$$
{\mathbb P}(U=V)=-\int_0^1 (\theta'\phi^2)(t)dt.  
$$
Besides, the copula $C_{\theta,\phi}$ has no singular component
if and only if $\theta$ is a constant function.  This case
is described more precisely in Section~\ref{secsub}.

          \section{Measures of association}
\label{measure}

In the next two sections, we note $(X,Y)$ a random pair with joint cdf $H$, copula $C$ and margins $F$ and $G$. The case $C=C_{\theta,\phi}$ is explicitly precised. 

          \subsection{Spearman's Rho}

Several invariant to strictly increasing function measures of association between the components of the random pair $(X,Y)$
can be considered:
the normalized volume between graphs of $H$ and $FG$~\cite{SIGMA},
Kendall's Tau~\cite{NELSEN06}, paragraph~5.1.1,
Gini's coefficient~\cite{NELSEN06}, Blomqwist's medial
correlation coefficient~\cite{NELSEN06}, paragraph~5.1.4,
and Spearman's Rho~\cite{NELSEN06}, paragraph~5.1.2,
which is the probability of concordance minus 
the probability of discordance of two random pairs  
with joint cdf $H$ and $FG$.
Here, we focus on this latter measure, 
showing in Subsection~\ref{subsecfin} that this
measure can achieve any value in $[-3/4,1]$.
A first step towards this result consists in remarking
that Spearman's Rho can be written only in terms of the copula $C$:
$$
\rho=12\int_0^1 \!\!\! \int_0^1 C(u,v) dudv\!-\!3.
$$
Note that $\rho$ coincides with the correlation coefficient between the uniform
marginal distributions.
In the case of a copula generated by (\ref{defcop}), 
it can be expressed thanks to the functions $\phi$ and $\theta$.
\begin{prop}
\label{proprho}
Let $(X,Y)$ be a random pair with copula $C_{\theta,\phi}$.
The Spearman's Rho is given by
$$
\rho_{\theta,\phi}=12\left[\Phi ^2(1)\theta(1) -\int_0^1\Phi^2(t)\theta'(t)dt\right],
$$
where $\Phi(t)=\int_0^t \phi(x)dx$.
\end{prop}

\noindent {\bf Proof: }
Clearly, Spearman's Rho can be expanded as
$$
\rho_{\theta,\phi} =  12\left[\int_0^1\int_v^1\phi(u)\phi(v)\theta(u)dudv +\int_0^1\int_0^v\phi(u)\phi(v)\theta(v)dudv \right].
$$
An integration by parts shows that both terms are equal and thus,
\begin{eqnarray*}
\rho_{\theta,\phi}  &=& 24\int_0^1\left(\phi(v)\theta(v)\int_0^v \phi(u)du \right)dv\\
		   &=& 24\int_0^1 \theta(v)\phi(v)\Phi(v)dv\\
		   &=& 12\left[\Phi ^2(1)\theta(1) -\int_0^1\Phi^2(t)\theta'(t)dt\right],
\end{eqnarray*}
by a new integration by parts.
\CQFD
             \subsection{Tail dependence}
The upper tail dependence can be quantified
by the upper tail dependence coefficient~\cite{JOE}, paragraph~2.1.10,  defined as: 
$$
\lambda=\lim_{t\to 1^-} {\mathbb P}(F(X)>t | G(Y)>t).
$$
Again, this coefficient can be written only in terms of the copula $C$
(see~\cite{NELSEN06}, Theorem~5.4.2):
$$
\lambda=\lim_{u\rightarrow 1^-}\frac{\bar{C}(u,u)}{1-u},
$$
where $\bar{C}$ is the survival copula, {\it i.e.} $\bar{C}(u,v)=1-u-v+C(u,v)$.
In our family, the following simplified expression can be obtained:
\begin{prop}
\label{proplambda}
Let $(X,Y)$ be a random pair with copula $C_{\theta,\phi}$.
The upper tail dependence coefficient is:
$
\lambda_{\theta,\phi}=-\phi^2(1)\theta'(1).
$
\end{prop}

\noindent {\bf Proof:} Clearly, the upper tail dependence coefficient can be simplified as
$$
\lambda_{\theta,\phi} = \lim_{u\rightarrow 1^-}\frac{\phi^2(u)\theta(u)}{1-u}.
$$
Taking into account of (b) yields
$$
\lambda_{\theta,\phi} = -\lim_{u\rightarrow 1^-}\frac{\phi^2(1)\theta(1)-\phi^2(u)\theta(u)}{1-u}
		     = -(\phi^2\theta)'(1) = -\phi^2(1)\theta'(1),
$$
and the result is proved.\CQFD

\noindent Note that a coefficient measuring the lower tail dependence
can also be defined as, 
$$
\lim_{t\to 0^+} {\mathbb P}(F(X)\leq t | G(Y)\leq t).
$$
but it is always zero in the considered family.

   \section{Concepts of dependence}
\label{dependence}

In this section, for the sake of simplicity, we assume that $X$ any $Y$
are exchangeable. 
Several concepts of positive dependence have been introduced
and characterized in terms of copulas. $X$ and $Y$ are
\begin{itemize}
\item Positively Quadrant Dependent (PQD) if 
${\mathbb P}(X\leq x,Y\leq y)\geq {\mathbb P}(X\leq x){\mathbb P}(Y\leq y)$,
for all $(x,y)\in \mathbb{R}^2$
or equivalently
\begin{equation}
\label{PQD}
\forall (u,v)\in I^2,\;\; C(u,v)\geq uv.
\end{equation}
\item Left Tail Decreasing (LTD) if ${\mathbb P}(Y\leq y| X\leq x)$ is non-increasing in $x$ for all $y$, or equivalently, 
see Theorem~5.2.5 in~\cite{NELSEN06}, 
$u\to C_{\theta,\phi}(u,v)/u$ is non-increasing for all $v\in I$.

\item Right Tail Increasing (RTI) if ${\mathbb P}(Y>y|X>x)$ is non-decreasing in $x$ for all~$y$ or, equivalently,
$u\to (v-C_{\theta,\phi}(u,v))/(1-u)$ is non-increasing for all $v\in I$.

\item Left Corner Set Decreasing (LCSD) if ${\mathbb P}(X \leq x,Y \leq y|X \leq x',Y \leq y')$ is non-increasing in $x'$ and $y'$ for all $x$ and $y$, or equivalently,
see  Corollary~5.2.17 in~\cite{NELSEN06}, 
$C$ is a totally positive function of order 2,
{\it i.e.  } for all $(u_1,u_2,v_1,v_2)\in I^4$ such that $u_1\leq u_2$
and $v_1\leq v_2$, one has
\begin{equation}
\label{LCSD2}
C(u_1,v_1)C(u_2,v_2)-C(u_1,v_2)C(u_2,v_1)\geq 0.
\end{equation}
\item Right Corner Set Increasing (RCSI) if ${\mathbb P}(X>x,Y >y|X>x',Y>y')$ is non-decreasing in $x'$ and $y'$ for all $x$ and $y$, or equivalently,
the survival copula $\hat C$ associated to C is a totally positive function of order 2.
\end{itemize}
Concepts of negative dependence can be similarly defined.   
Recall that $\theta$ is supposed not to be the null function on $I$
and introduce
$$
v^*=\sup\{v\in I; \; \theta(v)\neq 0\}.
$$
The point $v^*$, which can be seen as the endpoint of $\theta$,  
plays a central role in the dependence properties of the copula $C_{\theta,\phi}$,
see Theorem~\ref{theodepcop1} below.

\begin{theo}
\label{theodepcop1}
Let $(X,Y)$ a random pair with copula $C_{\theta,\phi}$.
$X$ and $Y$ are 
\begin{itemize}
\item [(i)] PQD if and only if $\theta(u)\geq 0$ for all $u\in I$
and $\phi$ has a constant sign on $[0,v^*]$.
\item [(ii)] LTD if and only if $\theta(u)\geq 0$ for all 
$u\in I$ and $ u\to\phi(u)/u$ and $u\to\theta(u)\phi(u)/u$ are either both non-increasing or both non-decreasing on $[0,v^*]$.
\item [(iii)] RTI if and only if $\theta(u)\geq 0$ for all 
$u\in I$ and $ u\to\phi(u)/(1-u)$ 
and $u\to\theta(u)\phi(u)/(1-u)$ are either both non-increasing or both non-decreasing on $[0,v^*]$.
\item [(iv)] LCSD if and only if they are LTD.
\item [(v)]  RCSI if and only if they are RTI.
\end{itemize}
\end{theo}

\noindent {\bf Proof: }
(i): Condition~(\ref{PQD}) can be rewritten as
\begin{equation}
\label{PQD2}
\theta(\max(u,v))\phi(u)\phi(v)\geq 0,\;\; \forall (u,v)\in I^2. 
\end{equation}
Suppose first that $(X,Y)$ is PQD.  Considering $u=v$ in~(\ref{PQD2})
shows that $\theta(u)\phi^2(u)\geq 0$ for all $u\in I$.  Since $\phi$
vanishes at most on isolated points, $\theta(u)\geq 0$ almost everywhere
on $I$.  Recalling that $\theta$ is continuous on $I$, we have
$\theta(u)\geq 0$ for all $u\in I$. Moreover, from (d),
$\theta$ in non-increasing on $I$, and consequently $\theta(t)>0$
for all $t\in [0,v^*)$. Thus, for all $(u,v)\in [0,v^*)^2$,
$\theta(\max(u,v))>0$ and
condition~(\ref{PQD2}) yields 
$\phi(u)\phi(v)\geq 0$ which implies that $\phi$
has a constant sign on $[0,v^*]$. \\ 
Conversely, suppose $\theta(u)\geq 0$ for all $u\in I$
and $\phi$ has a constant sign on $[0,v^*]$. For symmetry reasons,
it suffices to verify condition~(\ref{PQD2}) for $0\leq u\leq v\leq 1$.  
In this case $\theta(v)\phi(u)\phi(v)=0$ if $v\geq v^*$ and
$\theta(v)\phi(u)\phi(v)\geq 0$ otherwise.  

(ii) and (iii): Proofs are similar. Focusing on (iii), the necessary and sufficient condition can be rewritten as
\begin{equation}
\label{RTI2}
u\to \theta(\max(u,v))\phi(v)\phi(u)/(1-u) \mbox{ is non-decreasing for all } 
v\in [0,v^*].
\end{equation}
Supposing that $(X,Y)$ is RTI also implies that $(X,Y)$ is PQD
and consequently $\theta(u)\geq 0$ for all $u\in I$
and $\phi$ has a constant sign on $[0,v^*]$. Assuming for example
that $\phi(u)\geq 0$ for all $u\in[0,v^*]$,  condition~(\ref{RTI2})
implies that $u\to\phi(u)/(1-u)$ is non-decreasing on $[0,v]$ and
that $u\to\theta(u)\phi(u)/(1-u)$ is non-decreasing on $[v,v^*]$,
for all $v\in[0,v^*]$.  Thus, $u\to\phi(u)/(1-u)$ and 
$u\to\theta(u)\phi(u)/(1-u)$ are both non-decreasing on $[0,v^*]$.  
Conversely, assume $\theta(u)\geq 0$ for all $u\in I$, and
that $u\to\phi(u)/(1-u)$ and $u\to\theta(u)\phi(u)/(1-u)$ are
non-decreasing on  $[0,v^*]$.  
From Lemma~\ref{lemme0}(i) in the appendix, $\phi$ is non-negative on $[0,v^*]$,  
and~(\ref{RTI2}) is clearly true.  

(iv) and (v): Proofs are similar. Let us focus on (iv). It is well-known that $(X,Y)$ LCSD implies $(X,Y)$ LTD.
Let us prove that the converse result is also true in the $C_{\theta,\phi}$ family.
Suppose $(X,Y)$ is LTD. Following~(ii), one can assume that
$u\to\phi(u)/u$ is non-increasing and $\phi(u)\geq 0$
for all $u\in[0,v^*]$, together with $\theta(u)\geq 0$ for all $u\in I$.  
Lemma~\ref{lemme0}(ii) entails that $u\to\theta(u)\phi(u)/u$ is non-increasing 
on $[0,v^*]$.  
Four cases have to be considered to prove~(\ref{LCSD2}):\\
-- If $0\leq u_1 \leq u_2\leq v_1\leq v_2\leq 1$, condition~(\ref{LCSD2})
reduces to $A_1\geq 0$, where we have defined
$$
A_1:=\left[\frac{\theta(v_1)\phi(v_1)}{v_1}-\frac{\theta(v_2)\phi(v_2)}{v_2}\right]\left[\frac{\phi(u_1)}{u_1}-\frac{\phi(u_2)}{u_2}\right].
$$
If $u_2\leq v^*$, Lemma~\ref{lemme00} in the appendix yields $A_1\geq 0$.
Otherwise, $u_2> v^*$ implies $v_1\geq v^*$ and Lemma~\ref{lemme00} yields $A_1= 0$.\\
-- If $0\leq u_1 \leq v_1\leq u_2\leq v_2\leq 1$, condition~(\ref{LCSD2})
can be rewritten $A_2\geq 0$, with
$$
A_2:=A_1+C_{\theta,\phi}(u_1,v_2)[\theta(v_1)-\theta(u_2)] 
\frac{\phi(u_2)\phi(v_1)}{u_1 u_2 v_1 v_2}.
$$
If $u_2\leq v^*$, then $\phi(u_2)\phi(v_1)\geq 0$ and Lemma~\ref{lemme00} yields $A_1\geq 0$.
Consequently, $A_2\geq 0$.\\
If $v_1\leq v^*\leq u_2$, then $\theta(u_2)=\theta(v_2)=0$ and $A_2$ reduces to
$$
A_2= \frac{\theta(v_1)\phi(v_1)\phi(u_1)}{u_1 v_1 }\geq 0.
$$
Finally, if $v^*\leq v_1$, then $A_2=0$.\\
-- If $0\leq u_1 \leq v_1\leq v_2\leq u_2\leq 1$, condition~(\ref{LCSD2})
can be rewritten $A_3\geq 0$, with
\begin{eqnarray*}
A_3:= A_1 + \frac{\phi(u_2)}{u_2}&& \left[ 
\frac{\phi(v_1)}{v_1} \frac{C_{\theta,\phi}(u_1,v_2)}{u_1v_2}(\theta(v_1)-\theta(u_2))
+
\frac{\phi(v_2)}{v_2} \frac{C_{\theta,\phi}(u_1,v_1)}{u_1v_1}(\theta(u_2)-\theta(v_2))
\right]\\
=A_1 + \frac{\phi(u_2)}{u_2}&& \left[
 (\theta(v_1)-\theta(u_2))  \left(
\frac{\phi(v_1)}{v_1} - \frac{\phi(v_2)}{v_2} \right)\right.\\
&&+ \left.(\theta(v_1)-\theta(v_2)) \frac{\phi(v_2)}{v_2}
\left( 1 +\theta(u_2) \frac{\phi(u_1)}{u_1}\frac{\phi(v_1)}{v_1}\right)\right].
\end{eqnarray*}
If $u_2\leq v^*$, then Lemma~\ref{lemme00} yields $A_1\geq 0$ and all the above
differences are non-negative.
Consequently, $A_3\geq 0$.\\
If $v_2\leq v^*\leq u_2$, then $\theta(u_2)=0$ and $A_3$ reduces to
$$
A_3= \left(\frac{\theta(v_1)\phi(v_1)}{v_1}- \frac{\theta(v_2)\phi(v_2)}{v_2}\right)
\frac{\phi(u_1)}{u_1} \geq 0.
$$
If $v_1\leq v^*\leq v_2$, then $\theta(u_2)=\theta(v_2)=0$ and $A_3$ reduces to
$$
A_3= \frac{\theta(v_1)\phi(v_1)}{v_1}
\frac{\phi(u_1)}{u_1} \geq 0.
$$
Finally, if $v^*\leq v_1$, then $A_3=0$.\\
-- The three remaining situations 
are equivalent to the three previous ones
since the considered copulas are symmetric in the arguments.
\CQFD

\section{Sub-families and examples}
\label{secsub}

Recall that (b) is true if $\phi(1)=0$ 
or $\theta(1)=0$. The corresponding
sub-families are now studied in details and examples of copulas
in each sub-family are given.  

   \subsection{The case $\theta(1)=0$}

Let us focus on the sub-family of $C_{\theta,\phi}$ defined by conditions
(a), (b1), (c) and (d), where
\begin{itemize}
\item [(b1)] $\theta(1)=0$,
\end{itemize}
First, note that (b1, d) implies that $\theta$ is non negative on $I$.
From Proposition~\ref{proprho} and (b1), Spearman's Rho is
given by
\begin{equation}
\label{rhosub1}
\rho_{\theta,\phi}= -12\int_0^1\Phi^2(t)\theta'(t)dt,
\end{equation}
and (d) entails that, in this sub-family, $\rho_{\theta,\phi}\geq 0$.  
Second, we focus on copulas generated by univariate cdf 
and defined by
$$
C_{{\bar{K}}^{-1},Id}(u,v)=uv [ 1+ {\bar{K}}^{-1}(\max(u,v))],
$$
where $K$ is a cdf on ${\mathbb R}^+$,
$\bar{K}$ is the associated survival function,
$\bar{K}^{-1}$ is its generalized inverse defined as
$\bar{K}^{-1}(x)=K^{-1}(1-x)=\inf\{t\geq 0,\; K(t)\geq 1-x\}$
and $\phi=Id$ is the identity function.
We assume that $K$ 
is strictly increasing and differentiable on $(K^{-1}(0),K^{-1}(1))$,
the associated point distribution function is denoted by $k$.
The following corollary provides sufficient and necessary conditions
to ensure that $C_{{\bar{K}}^{-1},Id}$ is a copula.
It shows that the hazard function $k/\bar K$ is the key quantity
in this context.
\begin{coro}
\label{corofin}
$C_{{\bar{K}}^{-1},Id}$ is a copula if and only if, 
for all $t\geq 0$ such that $0<K(t)<1$,
\begin{equation}
\label{eqtail}
\frac{k(t)}{\bar{K}(t)}\geq \frac{1}{1+t}.
\end{equation}
\end{coro}
\noindent {\bf Proof: } Condition (a) is verified
since $\phi(x)=x$. Besides, $K(0)=0$ is equivalent to
condition (b). Condition (c) is
equivalent to 
$$
\frac{x}{k\left(\bar{K}^{-1}(x)\right)}-\bar{K}^{-1}(x)\leq 1,
$$
for all $x\in I$, which can be rewritten as 
$$
\frac{\bar{K}(t)}{k(t)}-t\leq 1,
$$
for all $t\geq 0$ by introducing $t=\bar{K}^{-1}(x)$.
The conclusion follows from Theorem~\ref{theocop}.
\CQFD
\noindent
As a consequence of condition~(\ref{eqtail}), one can easily
show that necessarily, $\bar{K}(x)\leq 1/(1+x)$.
Let us also note that, from~(\ref{rhosub1}) and Proposition~\ref{proplambda},
Spearman's Rho as well as the upper tail dependence coefficient can be rewritten in terms of $K$ as
\begin{equation}
\label{rhosub2}
\rho_{\bar{K}^{-1},Id}= 3\int_0^{+\infty} \bar{K}^4(t)dt\; \mbox{ and }\;
\lambda_{\bar{K}^{-1},Id}=1/k(0).
\end{equation}
In this sub-family, Blomqwvist's medial correlation coefficient $\beta$
benefits of a nice interpretation
$$
\beta_{\bar{K}^{-1},Id}=4C_{\bar{K}^{-1},Id}(1/2,1/2)-1={\bar K}^{-1}(1/2)=K^{-1}(1/2),
$$
as the median of the cdf $K$. 
Besides, characterizations of dependence properties in Proposition~\ref{theodepcop1} can be simplified as
\begin{coro}
Let $(X,Y)$ a random pair with copula $C_{\bar{K}^{-1},Id}$.
$X$ and $Y$ are always PQD, LTD and LCSD.
Moreover, $X$ and $Y$ are RTI and RCSI if and only if
for all $t>0$ such that $0<K(t)<1$,
$$
\frac{k(t)}{K(t)\bar{K}(t)}\geq \frac{1}{t}.
$$
\end{coro}
The proof is similar to the one of Corollary~\ref{corofin}.
In examples~1,...4, all the PQD, LTD, LCSD, RTI and RCSI
properties hold.

\noindent {\bf Examples}

1. A first example of copula belonging to this sub-family is
the Cuadras-Aug\'e copula~\cite{CA81}:  
$$
C_\alpha^{\mbox{\tiny CA}}(u,v)=\min(u,v)^{\alpha}(uv)^{1-\alpha}=M^\alpha(u,v) \Pi^{1-\alpha}
(u,v),
$$
where $\alpha\in[0,1]$, $M$ is the Fr\'echet upper bound defined by
$M(u,v)=\min(u,v)$ and $\Pi$ is the product copula $\Pi(u,v)=uv$.  
The copula $C_\alpha^{\mbox{\tiny CA}}$ can be interpreted as the weighted geometric mean
of $M$ and $\Pi$. It is generated by the $C_{\bar{K}^{-1},Id}$ family
with $\bar{K}(x)=(1+x)^{-1/\alpha}$, which is the survival function of a 
Generalized Pareto Distribution (GPD)
with positive shape parameter $1/\alpha$ (see for instance Table~1.2.6 in~\cite{EMBR}). The
associated Spearman's Rho given by~(\ref{rhosub2}) is 
$\rho_\alpha^{\mbox{\tiny CA}}=3\alpha/(4-\alpha)$ and the upper tail dependence coefficient
is $\lambda^{\mbox{\tiny CA}}_\alpha=\alpha$.

2. Another similar example is the family (B11), introduced in~\cite{JOE}, page~148:
$$
C_\sigma^{\mbox{\tiny B11}}(u,v)=\sigma\min(u,v)+(1-\sigma)uv=\sigma M(u,v) + (1-\sigma)\Pi (u,v),
$$
where $\sigma\in(0,1]$. The copula $C_\sigma^{\mbox{\tiny B11}}$ can be interpreted as the
 weighted arithmetic mean of $M$ and $\Pi$. 
It is generated by the $C_{\bar{K}^{-1},Id}$ family
with $\bar{K}(x)=(1+ x/\sigma)^{-1}$, which is the survival function of a GPD
with scale parameter $\sigma$ (see for instance Table~1.2.6 in~\cite{EMBR}).
The associated Spearman's Rho 
and upper tail dependence coefficient are
$\rho_\sigma^{\mbox{\tiny B11}}=\lambda^{\mbox{\tiny B11}}_\sigma=\sigma$.
Note that, for all $\alpha\in(0,1]$, one always has 
$\rho_\alpha^{\mbox{\tiny B11}}\geq \rho_\alpha^{\mbox{\tiny CA}}$.
Since both Cuadras-Aug\'e and (B11) copulas are indexed by a single
parameter they do not allow the pair $(\rho,\lambda)$ to reach arbitrary values 
in $[0,1]^2$. 

3. To partially overcome this limitation, it is natural to consider
$\bar{K}(x)=(1+x/\sigma)^{-1/\alpha}$, which is the survival function
of a GPD with positive
shape parameter $1/\alpha$ and scale parameter
$\sigma$, $\alpha\in(0,1]$, $\alpha\sigma\in(0,1]$. The
associated Spearman's Rho given by~(\ref{rhosub2}) is 
$\rho_{\alpha,\sigma}^{\mbox{\tiny GPD}}=3\alpha\sigma/(4-\alpha)$
 and the upper tail dependence coefficient
 $\lambda^{\mbox{\tiny GPD}}_{\alpha,\sigma}=\alpha\sigma$.
Thus, the $C_{\alpha,\beta}^{\mbox{\tiny GPD}}$ copula allows the pair
$(\rho,\lambda)$ to reach any value in the triangle
$\{ (\rho,\lambda)\in(0,1)^2\; :\; \rho\leq\lambda<4\rho/3\}$
with the following choice of parameters: $\alpha(\lambda,\rho)=4-3\lambda/\rho$
and $\sigma(\lambda,\rho)=(\rho\lambda)/(4\rho-3\lambda)$.

4. Choosing $K$ as the cdf
of the uniform distribution on $[0,\alpha]$, $\alpha\leq 1$
gives rise to the family of copulas
$$
C_\alpha^{\mbox{\tiny Uniform}}(u,v)=uv(1+\alpha\min(1-u,1-v)),
$$
introduced in~\cite{HUE}, Section~1, and with associated Spearman's Rho
$\rho_{\alpha}^{\mbox{\tiny Uniform}}=3\alpha/5$
and  upper tail dependence coefficient
 $\lambda^{\mbox{\tiny Uniform}}_{\alpha}=\alpha$.

5. Finally, note that the family
$$
C_f(u,v)=\min(u,v) f(\max(u,v))
$$
proposed in~\cite{DUR} can also enter our sub-family with
an appropriate choice of $K$.

\noindent Basing on Example~1, we can state the following result:
\begin{prop}
\label{proprhopos}
Suppose $C_{\theta,\phi}$ is a copula and $\theta(1)=0$.  
Thus, $0\leq \rho_{\theta,\phi}\leq 1$ and $0\leq \lambda_{\theta,\phi}\leq 1$,
and these bounds are reached within the sub-family.  
\end{prop}

           \subsection{The case $\phi(1)=0$}

Here, we focus on the sub-family of $C_{\theta,\phi}$ defined by conditions
(a), (b2), (c) and (d), where
\begin{itemize}
\item [(b2)] $\phi(1)=0$,
\end{itemize}
Note that (b2) implies that the upper tail dependence coefficient is
always  null in this sub-family.  
This sub-family encompasses the semiparametric family of copulas
with constant function $\theta$ defined in~(\ref{eqcopunous}).
Consequently, this sub-family also
includes the FGM family~(\ref{eqFGM}),
the parametric family 
of symmetric copulas with cubic sections proposed in~\cite{NELSEN97},
equation~(4.4), both kernel families~(\ref{eqHK1}) and~(\ref{eqHK2})
introduced in~\cite{HUE},
and the PQD copulas~(\ref{eqLX}) introduced in~\cite{LAI}.\\
From Proposition~\ref{proprho}, in the subfamily
of $C_{\theta,\phi}$ constrained by (b2), the following
lower bound for Spearman's Rho holds:
$$
\rho_{\theta,\phi}\geq 12\Phi^2(1)\theta(1),
$$
where the right-hand term can be interpreted as Spearman's Rho 
associated to the copula~(\ref{eqcopunous}) with constant function
$\theta(.)=\theta(1)$. Since, in this particular case, Spearman's 
Rho is lower bounded by $-3/4$ (see~\cite{Nous}, Proposition~2),
we have:
\begin{prop}
\label{propfin}
Suppose $C_{\theta,\phi}$ is a copula and $\phi(1)=0$.  
Thus, $\lambda_{\theta,\phi}= 0$ and $\rho_{\theta,\phi}\geq -3/4$,
and this bound is reached within the subfamily.  
\end{prop}

\noindent {\bf Remark.}
It is of course possible to build copulas such that $\phi(1)=0$
and $\theta$ is a non constant function. As an example, consider
the function $\phi(x)=x(1-x)$ which generates the FGM family of copulas.
Taking $u=0$ in condition~(c) and integrating with respect to $v\in[x,1]$
imply that $\theta(x)\leq 1/x$ for all 
$0<x\leq 1$. Let us consider the extreme case $\theta(x)=1/x$.
The copula writes
$$
C(u,v)=\Pi(u,v) + (1-u)(1-v)M(u,v),
$$
and the associated Spearman's Rho is $\rho=3/5$ which is much
larger than the maximum value ($\rho=1/3)$ in the FGM family.

\subsection{General case}
\label{subsecfin}

Collecting Proposition~\ref{proprhopos} and Proposition~\ref{propfin},
we are now in position to provide the bounds for the
general family~(\ref{defcop}).
\begin{prop}
Suppose $C_{\theta,\phi}$ is a copula.
Thus, $0\leq \lambda_{\theta,\phi} \leq 1$ and $-3/4\leq \rho_{\theta,\phi}\leq 1$, and these bounds are reached within the family.  
\end{prop}
Besides, Proposition~\ref{propsing} entails
that the copulas~(\ref{eqcopunous}) are the only ones which are 
absolutely continuous.  Thus, we can conclude that, in the general
$C_{\theta,\phi}$ family, the absolute continuity is incompatible
with the upper tail dependence.


\section*{Appendix: Auxiliary lemmas}

\begin{lemm}
\label{lemme0}
Assume $C_{\theta,\phi}$ is a copula.
\begin{itemize}
\item [(i)] If $u\to\phi(u)/(1-u)$ is non-decreasing (resp. non-increasing)
on $[0,v^*]$ then $\phi(u) \geq 0$ (resp. $\leq 0$) for all $u\in [0,v^*]$.
\end{itemize}
If, moreover, $\theta(u)\geq 0$ for all $u\in I$, then
\begin{itemize}
\item [(ii)] If $\phi(u)\leq 0$ (resp. $\geq 0$) for all $u\in J\subset I$ and
$u\to\phi(u)/u$ is non-decreasing (resp. non-increasing) on $J$ then
$u\to\theta(u)\phi(u)/u$ is non-decreasing (resp. non-increasing) on
$J$.
\end{itemize}
\end{lemm}
\noindent {\bf Proof:}
\begin{itemize}
\item [(i)] Assume $u\to\phi(u)/(1-u)$ is non-decreasing on $[0,v^*]$. Then,
from (a), $\forall u \in [0,v^*]$, $\phi(u)/(1-u)\geq \phi(0)=0$.
Therefore, $\phi(u)/(1-u)$ is non-negative on $[0,v^*]$ and the conclusion
follows.
\item [(ii)] Remark that
$[\theta(u)\phi(u)/u]'=\theta'(u)\phi(u)/u+\theta(u)[\phi(u)/u]'$. Thus, if
$\phi(u)\leq 0$, $\theta(u)\geq 0$ and $u\to \phi(u)/u$ is non-decreasing,
(d) implies that $[\theta(u)\phi(u)/u]'\geq 0$ for all $u\in J$.
\CQFD
\end{itemize}

\begin{lemm}
\label{lemme00}
Assume $C_{\theta,\phi}$ is a copula, $\theta(u)\geq 0$ for all
$u\in I$ and
\begin{itemize}
\item [(i)] either $\{ u\to\phi(u)/u$ is non-increasing and $\phi(u)\geq 0 \}$
for all $u\in[0,v^*]$,
\item [(ii)] or $\{u\to\phi(u)/u$ is non-decreasing and $\phi(u)\leq 0\}$
for all $u\in[0,v^*]$.
\end{itemize}
Let $(u_1,u_2,v_1,v_2)\in I^4$ such that $u_1\leq u_2$, $v_1\leq v_2$
and introduce 
$$
A_1:=\left[\frac{\theta(v_1)\phi(v_1)}{v_1}-\frac{\theta(v_2)\phi(v_2)}{v_2}\right]\left[\frac{\phi(u_1)}{u_1}-\frac{\phi(u_2)}{u_2}\right].
$$
Then, $u_2\leq v^*$ entails $A_1\geq 0$ and $v^*\leq v_1$ entails $A_1=0$.
\end{lemm}
\noindent {\bf Proof:} Assume (i) holds, situation (ii) is similar.
First, remark that $u\to\theta(u)\phi(u)/u$ is non-increasing on the whole $I$
interval, since this function is 
non-negative and 
non-increasing on $[0,v^*]$ (Lemma~\ref{lemme0}(ii))
and vanishes on $[v^*,1]$. Therefore,
$$
\frac{\theta(v_1)\phi(v_1)}{v_1}-\frac{\theta(v_2)\phi(v_2)}{v_2}\geq 0
$$
in all cases. Now, if $u_2\leq v^*$, then $u\to\phi(u)/u$ is non-increasing
on the considered interval and the conclusion follows.
It $v^*\leq v_1$, then $\theta(v_1)=\theta(v_2)=0$ and the result is proved.
\CQFD
\end{document}